\documentclass[12pt]{article}
\usepackage{amsmath}
\usepackage{amssymb}
\usepackage{amscd}
\usepackage{amsthm}

\theoremstyle{plain}
\newtheorem{Thm}{Theorem}[section]
\newtheorem{Conj}[Thm]{Conjecture}

\newtheorem{Lem}[Thm]{Lemma}

\theoremstyle{definition}
\newtheorem{Defn}[Thm]{Definition}

\theoremstyle{remark}
\newtheorem{Rem}[Thm]{Remark}

\numberwithin{equation}{section}

\title{On a Relative Version of Fujita's Freeness Conjecture}
\author{Yujiro Kawamata}
\date{\today}

\begin{document}

\maketitle

\section{Introduction}

The following is Fujita's freeness conjecture:

\begin{Conj}\label{absolute1}
Let $X$ be a smooth projective variety of dimension $n$ and
$H$ an ample divisor.
Then the invertible sheaf $\mathcal{O}_X(K_X + mH)$
is generated by global sections if $m \ge n+1$, or $m = n$ and $H^n \ge 2$.
\end{Conj}

We have a stronger local version of Conjecture~\ref{absolute1}
(cf. \cite{fujita}):

\begin{Conj}\label{absolute2}
Let $X$ be a smooth projective variety of dimension $n$,
$L$ a nef and big invertible sheaf on $X$, and $x \in X$ a point.
Assume that
$L^n > n^n$ and $L^dZ \ge n^d$ for any irreducible subvariety $Z$ of $X$
of dimension $d$ which contains $x$.
Then the natural homomorphism
\[
H^0(X, \omega_X \otimes L) \to \omega_X \otimes L \otimes \kappa(x)
\]
is surjective.
\end{Conj}

We shall extend the above conjecture to a relative setting.

Let $f: Y \to X$ be a surjective morphism of smooth projective varieties.
We note that the geometric fibers of $f$ are not necessarily connected.
Assume that there exists a normal crossing divisor $B = \sum_{i=1}^h B_i$ 
on $X$ 
such that $f$ is smooth over $X_0 = X \setminus B$.
Then the sheaves $R^qf_*\omega_{Y/X}$ are locally free for $q \ge 0$
(\cite{abel} for $q = 0$ and \cite{Kollar} in general).
We note that even if we change the birational model of $Y$,
the sheaf $R^qf_*\omega_{Y/X}$ does not change.

The relative version is the following:

\begin{Conj}\label{relative1}
Let $f: Y \to X$ be a surjective morphism from a smooth projective variety to 
a smooth projective variety of dimension $n$ 
such that $f$ is smooth over $X_0 = X \setminus B$ for 
a normal crossing divisor $B$ on $X$.
Let $H$ be an ample divisor on $X$.
Then the locally free sheaf $R^qf_*\omega_Y \otimes \mathcal{O}_X(mH)$
is generated by global sections if $m \ge n+1$, or $m = n$ and $H^n \ge 2$.
\end{Conj}

We have again a stronger local version:

\begin{Conj}\label{relative2}
Let $f: Y \to X$ be a surjective morphism from a smooth projective variety to 
a smooth projective variety of dimension $n$
such that $f$ is smooth over $X_0 = X \setminus B$ for 
a normal crossing divisor $B$ on $X$.
Let $L$ be a nef and big invertible sheaf on $X$, and $x \in X$ a point.
Assume that $L^n > n^n$ and $L^dZ \ge n^d$ for any irreducible subvariety $Z$ of $X$
of dimension $d$ which contains $x$.
Then the natural homomorphism
\[
H^0(X, R^qf_*\omega_Y \otimes L) \to R^qf_*\omega_Y \otimes L \otimes \kappa(x)
\]
is surjective for any $q \ge 0$.
\end{Conj}

In \cite{fujita}, the following strategy toward Conjectures~\ref{absolute1}
and \ref{absolute2} was developed:

\begin{Thm}[\cite{fujita}]\label{method}
Let $X$ be a smooth projective variety of dimension $n$,
$L$ a nef and big invertible sheaf, and $x \in X$ a point.
Assume the followig condition:
for any effective $\mathbb{Q}$-divisor $D_0$ on $X$ 
such that $(X, D_0)$ is KLT, 
there exists an effective $\mathbb{Q}$-divisor $D$ on $X$ 
such that 

(1) $D \equiv \lambda L$ for some $0 < \lambda < 1$, 

(2) The pair $(X, D_0 + D)$ is properly log canonical at $x$, and

(3) $\{x\}$ is a log canonical center for $(X, D_0 + D)$.

\noindent
Then the natural homomorphism
\[
H^0(X, \omega_X \otimes L) \to \omega_X \otimes L \otimes \kappa(x)
\]
is surjective.
\end{Thm}

\begin{Thm}[\cite{fujita}]\label{lowdim}
(1) In the situation of Conjecture~\ref{absolute2},  assume that $\dim X \le 3$.
Then the condition of Theorem~\ref{method} is satisfied
for $L$.

(2) In the situation of Conjecture~\ref{absolute1}, assume that $\dim X = 4$.
Then the condition of Theorem~\ref{method} is satisfied
for $L = mH$ with $m \ge 5$ or $m = 4$ and $H^4 \ge 2$. 
\end{Thm}

Our main result is the relative version of Theorem~\ref{method}:

\begin{Thm}\label{main}
Let $f: Y \to X$ be a surjective morphism from a smooth projective variety to 
a smooth projective variety of dimension $n$
such that $f$ is smooth over $X_0 = X \setminus B$ for 
a normal crossing divisor $B$ on $X$.
Let $L$ be a nef and big invertible sheaf on $X$, and $x \in X$ a point.
Assume the followig condition:
for any effective $\mathbb{Q}$-divisor $D_0$ on $X$ 
such that $(X, D_0)$ is KLT, 
there exists an effective $\mathbb{Q}$-divisor $D$ on $X$ 
such that 

(1) $D \equiv \lambda L$ for some $0 < \lambda < 1$, 

(2) The pair $(X, D_0 + D)$ is properly log canonical at $x$, and

(3) $\{x\}$ is a log canonical center for $(X, D_0 + D)$.

\noindent
Then the natural homomorphism
\[
H^0(X, R^qf_*\omega_Y \otimes L) \to R^qf_*\omega_Y \otimes L \otimes \kappa(x)
\]
is surjective for any $q \ge 0$.
In particular, Conjecture~\ref{relative2} for $\dim X \le 3$ and
Conjecture~\ref{relative1} for $\dim X = 4$ hold.
\end{Thm}

The main tool for the proof of Theorem~\ref{method} was the $\mathbb{Q}$-divisor version 
of the Kodaira vanishing theorem, so-called Kawamata-Viehweg vanishing theorem.
In order to prove Theorem~\ref{main}, we shall extend the Koll\'ar vanishing theorem,
the relative version of the Kodaira vanishing theorem, to the $\mathbb{Q}$-divisor version.
We had a similar result already in \cite{mK}~Theorem 3.3, but we need more precise 
version.

In \S 2, we review the construction of the canonical extension for a variation of 
Hodge structures.
In order to describe the behavior of the Hodge bundles at infinity, we
shall introduce the notion of parabolic structures over arbitrary dimensional base 
in \S 3.
In \S 4, 
we shall prove the base change theorem of the parabolic structures (Lemma~\ref{BC})
and derive the correct $\mathbb{Q}$-divisor version of the Koll\'ar vanishing theorem
(Theorem~\ref{vanishing}).
The main result will be proved in \S 5.

\begin{Rem}
(1) It is easy to prove Conjecture~\ref{absolute1}
in the case $L$ is very ample.
\cite{Kollar} proved Conjecture~\ref{relative1} in the case $L$ is 
very ample.

(2) Let $\pi: P = \mathbb{P}(\mathcal{F}) \to X$ be the projective space bundle
associated to our locally free sheaf $\mathcal{F} = R^qf_*\omega_{Y/X}$, 
and $\mathcal{O}_P(H)$ the tautological invertible sheaf.
Then the sheaf $\mathcal{F} \otimes \mathcal{O}_X(K_X + L)$ is generated by global sections, if and only if $\mathcal{O}_P(H + \pi^*(K_X + L))$ is so.
Since $K_P = - rH + \pi^*(K_X + \mathrm{det}(\mathcal{F}))$ if
$r = \mathrm{rank}(\mathcal{F})$, Conjecture~\ref{relative1} is related to, 
but different from, Conjecture~\ref{absolute1}.  
\end{Rem}

\section{Review on the Hodge bundles} 

Let $X$ be a smooth projective variety, and $B$ a normal crossing divisor.
Let $H_{\mathbb{Z}}$ 
be a local system on $X_0 = X \setminus B$.
A {\em variation of Hodge structures} on $X_0$ 
is defined as a decreasing filtration
$\{\mathcal{F}^p_0\}$ of locally free subsheaves on 
$\mathcal{H}_0 = H_{\mathbb{Z}} \otimes \mathcal{O}_{X_0}$
which satisfy certain axioms (\cite{G}).
Assume in addition that all the local monodromies of $H_{\mathbb{Z}}$ 
around the branches of $B$ are unipotent.
Then we define a locally free sheaf $\mathcal{H}$ on $X$ called 
the {\em canonical extension} of $\mathcal{H}_0$ as follows.
Let $\{z_1, \cdots, z_n\}$ be local coordinates at a given point
$x \in X$ such that $B$ is defined by an equation $z_1 \cdots z_r$ near $x$.
Let $T_i$ be the monodromies of $H_{\mathbb{Z}}$ around the branches of $B$
defined by $z_i = 0$.
Let $v$ be a multi-valued flat section of $H_{\mathbb{Z}}$.
Then the expression
\[
s = \exp(- \frac 1{2\pi\sqrt{-1}}\sum_{i=1}^r \log T_i \log z_i)v
\]
is single-valued and gives a holomorphic section of $\mathcal{H}_0$,
where $\log T_i$ is defined by a finite power series of $T_i - 1$.
The canonical extension 
$\mathcal{H}$ is defined as a locally free sheaf on $X$ whose local basis
near $x$ consists of the $s$ when $v$ making a basis of $H_{\mathbb{Z}}$.
The filtration $\{\mathcal{F}^p_0\}$ extends to a filtration $\{\mathcal{F}^p\}$ 
of locally free subsheaves on $\mathcal{H}$
such that the quotients $\mathcal{H}/\mathcal{F}^p$ are also 
locally free for all $p$ (\cite{Schmid}).

Let $f: Y \to X$ be a surjective morphism from another 
smooth projective variety which is smooth over $X_0$.
Let $d = \dim Y - \dim X$, $Y_0 = f^{-1}(X_0)$ and $f_0 = f \vert_{Y_0}$.
Then $H_{\mathbb{Z}} = R^{d+q}f_{0*}\mathbb{Z}_{Y_0}$ is a 
variation of Hodge structures on $X_0$ for any $q \ge 0$. 
We know that 
the canonical extension $\mathcal{F}^d$ in this case 
coincides with the direct image sheaf $R^qf_*\omega_{Y/X}$
(\cite{abel} and \cite{Kollar}).

The following lemma is obvious:

\begin{Lem}\label{canonical}
Let $X$ be a smooth projective variety, and $B$ a normal crossing divisor.
Let $H_{\mathbb{Z}}$ be a variation of Hodge structures on 
$X_0 = X \setminus B$
whose local monodromies around the branches of $B$ are unipotent, and
$\mathcal{H}$ the canonical extension of $\mathcal{H}_0 = 
H_{\mathbb{Z}} \otimes \mathcal{O}_{X_0}$ on $X$.
Let $\pi: X' \to X$ be a generically finite and surjective morphism 
from a smoth projective variety such that
$B' = (\pi^*B)_{\mathrm{red}}$ is a normal crossing divisor.
Let $H'_{\mathbb{Z}} = \pi^*H_{\mathbb{Z}}$ be the induced 
variation of Hodge structures on $X'_0 = X' \setminus B'$, 
and $\mathcal{H}'$ the canonical extension of $\mathcal{H}'_0 = 
H'_{\mathbb{Z}} \otimes \mathcal{O}_{X'_0}$ on $X'$.
Then $\mathcal{H}' = \pi^*\mathcal{H}$.
\end{Lem}

\section{Parabolic structure in several variables}

We generalize the notion of 
parabolic structures on vector bundles (\cite{MS}) over 
higher dimensional base space:

\begin{Defn}\label{para}
Let $f: Y \to X$ be a surjective morphism of smooth projective varieties.
Assume that there exists a normal crossing divisor $B = \sum_{i=1}^h B_i$ 
on $X$ 
such that $f$ is smooth over $X_0 = X \setminus B$.
Fixing a nonnegative integer $q$, 
we define a {\em parabolic structure} on the sheaf $\mathcal{F} = R^qf_*\omega_{Y/X}$.
It is a decreasing filtration of subsheaves
$F^{t_1,\dots,t_h} = F^{t_1,\dots,t_h}(\mathcal{F}) \subset \mathcal{F}$ 
with multi-indices $t = (t_1,\dots,t_h)$
($t_i \in \mathbb{R}_{\ge 0}$) defined by
\[
\begin{split}
\Gamma(U, F^{t_1,\dots,t_h}(\mathcal{F})) 
= \{&s \in \Gamma(U, \mathcal{F}) \,\vert \, 
(\prod_i z_i^{-t_i})s 
\textrm{ is } L^2 \\ &\textrm{ with respect to the Hodge metric} \},
\end{split}
\]
where $z_i$ is a local equation of the branch $B_i$ on an open subset 
$U \subset X$.
\end{Defn}

\begin{Lem}\label{parabolic}

(1) $F^t \supset F^{t'}$ for $t \le t'$, i.e., $t_i \le t'_i$ for all $i$.

(2) $F^{t_1, \dots, t_i+\epsilon, \dots, t_h} = 
F^{t_1, \dots, t_i, \dots, t_h}$ for $0 < \epsilon \ll 1$.

(3) $F^{t_1,\dots,t_i+1,\dots,t_h} = 
F^{t_1,\dots,t_i,\dots,t_h} \otimes \mathcal{O}_X(-B_i)$.

(4) Let $Y_0 = f^{-1}(X_0)$, $f_0 = f \vert_{Y_0}$ and $d = \dim Y - \dim X$.
If all the local monodromies of $R^{d+q}f_{0*}\mathbb{Z}_{Y_0}$ 
around the branches of $B$ are unipotent,
then $F^t = F^0$ for any $t = (t_1, \dots, t_h)$ 
with $0 \le t_i < 1$.
\end{Lem}

\begin{proof}
(1) through (3) are obvious.
(4) follows from the fact that the growth of the Hodge metric 
is logarithmic in this case (\cite{abel}).
\end{proof}

\begin{Rem}
For negative values of the $t_i$, we can also define $F^t$ as subsheaves of 
$\mathcal{F} \otimes \mathcal{O}_X(mB)$ for sufficiently large $m$
by using Lemma~\ref{parabolic}~(3).
We shall also write $F^{\sum_i t_iB_i}$ instead of $F^{t_1,\dots,t_h}$.
\end{Rem}

\begin{Defn}
For a local section $s \in \Gamma(U, \mathcal{F})$, we define
its {\em order of growth} along $B$ by 
\[
\mathrm{ord}(s) = \sum_i \mathrm{ord}_i(s)B_i 
=\inf \{\sum_i (1-t_i)B_i \,\vert \, s \in 
\Gamma(U, F^{t_1,\dots,t_h}(\mathcal{F}))\}.
\]
We note that $s \not\in \Gamma(U, F^{B - \mathrm{ord}(s)}(\mathcal{F}))$, and
\[
\Gamma(U, F^{t_1,\dots,t_h}(\mathcal{F})) 
= \{s \in \Gamma(U, \mathcal{F})
\,\vert \, \mathrm{ord}(s) + \sum_i t_iB_i < B\}.
\]
\end{Defn}

There is a nice local basis of the sheaf $\mathcal{F} = R^qf_*\omega_{Y/X}$:

\begin{Lem}\label{basis}
At any point $x \in X$, 
there exists an open neighborhood $U$ in the classical topology 
and a $\Gamma(U, \mathcal{O}_U)$-free
basis $\{s_1, \dots, s_k\}$ of $\Gamma(U, \mathcal{F})$ such that
\[
(\prod_i z_i^{\llcorner t_i + \mathrm{ord}_i(s_1) \lrcorner}) s_1,
\dots, (\prod_i z_i^{\llcorner t_i + \mathrm{ord}_i(s_k) \lrcorner}) s_k
\]
generates $\Gamma(U, F^{t_1,\dots,t_h}(\mathcal{F}))$ for any $t$,
where the $z_i$ are local equations of the $B_i$ on $U$.
In particular, 
the sheaf $F^{t_1,\dots,t_h}(\mathcal{F})$ is locally free
for any $t$.
\end{Lem}

\begin{proof}
We shall prove that the filtration $F^t$ is determined 
by the local monodromies of the cohomology sheaf $R^{d+q}f_{0*}\mathbb{Z}_{Y_0}$ 
around the branches of $B$ which are known to be quasi-unipotent, 
where $d = \dim Y - \dim X$.

We take an open neighborhood $U$ of $x \in X$ in the classical  
topology which is isomorphic to a 
polydisk with coordinates $\{z_1, \dots, z_n\}$ such that
$B \cap U$ is defined by $z_1 \cdots z_r = 0$.
To simplify the notation, we write $X$ instead of $U$.
There exists a finite surjective and Galois morphism $\pi: X' \to X$ 
from a smooth variety which is etale over $X_0$ such that, 
for the induced morphism $f': Y' \to X'$ from a 
desingularization $Y'$ of the fiber product $Y \times_X X'$, 
the local system 
$R^{d+q}f'_{0*}\mathbb{Z}_{Y'_0}$ has unipotent local monodromies around the
branches of $B' = \pi^{-1}(B)$, 
where we set $X'_0 = \pi^{-1}(X_0)$, $Y'_0 = f^{\prime -1}(X'_0)$
and $f'_0 = f' \vert_{Y'_0}$.
Let $\sigma: Y' \to Y$ be the induced morphism.

We may assume that $X'$ is isomorphic to a polydisk centered at
a point $x' = \pi^{-1}(x)$
with coordinates $\{z'_1, \dots, z'_n\}$, 
and the morphism $\pi: X' \to X$ is given by 
$\pi^*z_i = z_i^{\prime m_i}$ for some positive integers $m_i$,
where $m_i = 1$ for $i > r$.
The Galois group $G = \mathrm{Gal}(X'/X)$ is isomorphic to
$\prod_i \mathbb{Z}/(m_i)$.
Let $g_1, \dots, g_r$ be generators of $G$ such that
$g_i^*z'_j = \zeta_{m_i}^{\delta_{ij}}z'_j$ 
for some roots of unity $\zeta_{m_i}$ of order $m_i$.

The group $G$ acts on the sheaves
$R^qf'_*\omega_{Y'/X'}$ and $\omega_{X'}$ equivariantly 
such that the invariant part
$(\pi_*(R^qf'_*\omega_{Y'/X'} \otimes \omega_{X'}))^G$ is isomorphic to 
$R^qf_*\omega_{Y/X} \otimes \omega_X$, because 
$(\sigma_*\omega_{Y'})^G = \omega_Y$ and 
$R^p\sigma_*\omega_{Y'} = 0$ for $p > 0$.

The vector space $R^qf'_*\omega_{Y'/X'} \otimes \kappa(x')$ 
is decomposed into simultaneous eigenspaces with respect to the action of $G$.
Let $s_{x'}$ be a simultaneous eigenvector such that
$g_i^*s_{x'} = \zeta_{m_i}^{a_i}s_{x'}$ for some $a_i$ with $0 \le a_i < m_i$.
Let $\bar s'$ be a section of $R^qf'_*\omega_{Y'/X'}$ which extends $s_{x'}$.
Then the section 
\[
s' = \frac 1{\prod_i m_i}\sum_{i=1}^r\sum_{k_i=0}^{m_i-1} 
\frac{(\prod_i g_i^{k_i})^*\bar s'}{\prod_i \zeta_{m_i}^{a_ik_i}}
\]
satisfies that $s'(x') = s_{x'}$ and $g_i^*s' = \zeta_{m_i}^{a_i}s'$.

On the other hand, $dz'_1 \wedge \cdots \wedge dz'_n = 
(\prod_i m_i^{-1}z_i^{\prime 1-m_i}) dz_1 \wedge \cdots \wedge dz_n$ is a 
generating section of $\omega_{X'}$.
Therefore, $(\prod_i z_i^{\prime -a_i})s'$
descends to a section $s$ of $R^qf_*\omega_{Y/X}$.
If the $s_{x'}$ varies among a 
basis of $R^qf'_*\omega_{Y'/X'} \otimes \kappa(x')$, then the 
corresponding sections $s$ make a basis of the locally free sheaf
$R^qf_*\omega_{Y/X}$.
 
We have $\mathrm{ord}_i(s) = a_i/m_i$,
since the Hodge metric on the sheaf $R^qf'_*\omega_{Y'/X'}$ has
logarithmic growth along $B'$.
Therefore, the sections $(\prod_i z_i^{\llcorner t_i + a_i/m_i \lrcorner})s$
form a basis of a locally free sheaf
$F^{t_1,\dots,t_h}(R^qf_*\omega_{Y/X})$.
\end{proof}

\begin{Rem}
The Hodge metric and the flat metric on the canonical extension of
the variation of Hodge structures 
$R^{d+q}f_{0*}\mathbb{Z}_{Y_0} \otimes \mathcal{O}_{X_0}$ 
coincide when restricted to
the subsheaf $R^qf_*\omega_{Y/X}$.
Therefore, the statement that 
the Hodge metric on the canonical extension has logarithmic
growth is easily proved for the sheaf $R^qf_*\omega_{Y/X}$.
\end{Rem}

\section{Base change and a relative vanishing theorem}

By using the basis obtained in Lemma~\ref{basis}, we can study the 
base change property of the sheaf $R^qf_*\omega_{Y/X}$:

\begin{Lem}\label{BC}
Let $\pi: X' \to X$ be a generically finite and surjective morphism 
from a smooth projective variety such that
$B' = (\pi^*B)_{\mathrm{red}} 
= \sum_{i'=1}^{h'} B'_{i'}$ is a normal crossing divisor.
Let $\mu: Y' \to Y \times_X X'$ be a birational morphism from a smooth 
projective variety such that the induced morphism
$f': Y' \to X'$ is smooth over $X'_0 = X' \setminus B'$.
Let $\sigma: Y' \to Y$ be the induced morphism.
Then the following hold.

(1) Let $\{s_1, \dots, s_k\}$ be the basis of $\Gamma(U, R^qf_*\omega_{Y/X})$
in Lemma~\ref{basis}, and
let $U'$ be an open subset of $X'$ in the classical topology such that
$\pi(U') \subset U$.
Then the equality $\mathrm{ord}(\pi^*s_j) = \pi^*\mathrm{ord}(s_j)$ holds, and 
the basis 
$\{\pi^*s_1, \dots, \pi^*s_k\}$ of 
$\Gamma(U', \pi^*R^qf_*\omega_{Y/X})$
satisfies the conclusion of 
Lemma~\ref{basis} in the sense that sections 
\[
(\prod_{i'} z_{i'}^{\prime \llcorner t_{i'} + 
\mathrm{ord}_{i'}(\pi^*s_1) \lrcorner}) \pi^*s_1, \dots,
(\prod_{i'} z_{i'}^{\prime \llcorner t_{i'} + 
\mathrm{ord}_{i'}(\pi^*s_k) \lrcorner}) \pi^*s_k
\]
form a basis of 
$\Gamma(U', F^{t'_1,\dots,t'_{h'}}(R^qf'_*\omega_{Y'/X'}))$
for any $t' = (t'_1, \dots, t'_{h'})$,
where the $z'_{i'}$ are local equations of the $B'_{i'}$ on $U'$.
In particular, the sections 
\[
(\prod_{i'} z_{i'}^{\prime \llcorner 
\mathrm{ord}_{i'}(\pi^*s_1) \lrcorner}) \pi^*s_1, \dots,
(\prod_{i'} z_{i'}^{\prime \llcorner 
\mathrm{ord}_{i'}(\pi^*s_k) \lrcorner}) \pi^*s_k.
\]
form a basis of $\Gamma(U', R^qf'_*\omega_{Y'/X'})$.

(2) There is an equality of subsheaves of $\pi^*R^qf_*\omega_{Y/X}$:
\[
\begin{split}
&F^{t'_1,\dots, t'_{h'}}(R^qf'_*\omega_{Y'/X'}) \\
&= \sum_t \pi^*F^{t_1, \dots, t_h}(R^qf_*\omega_{Y/X}) 
\otimes \mathcal{O}_{X'}
(- \llcorner \sum_{i'} t'_{i'}B'_{i'} + \sum_i (1-t_i)\pi^*B_i \lrcorner).
\end{split}
\]
In particular,
\[
R^qf'_*\omega_{Y'/X'} = \sum_t \pi^*F^{t_1,\dots,t_h}(R^qf_*\omega_{Y/X}) 
\otimes \mathcal{O}_{X'}(- \llcorner \sum_i (1-t_i)\pi^*B_i \lrcorner).
\]
\end{Lem}

\begin{proof}
(1) Since the $s_j$ are derived from the basis in the case of unipotent 
monodromies, we obtain our assertion by Lemma~\ref{canonical}.

(2) We can check the assertion locally.
We write $\pi^*B_i = \sum_{i'} m_{ii'}B_{i'}$ for some nonnegative
integers $m_{ii'}$.
Then the left hand side is generated by the sections
\[
(\prod_{i'} z_{i'}^{\prime \llcorner t'_{i'} + 
\mathrm{ord}_{i'}(\pi^*s_j) \lrcorner}) \pi^*s_j
\]
for $1 \le j \le k$, while each component of the right hand side is by
\[
(\prod_{i'} z_{i'}^{\prime 
\sum_i (\llcorner t_i + \mathrm{ord}_i(s_j) \lrcorner)m_{ii'}
+ \llcorner t'_{i'} + \sum_i (1-t_i)m_{ii'} \lrcorner}) \pi^*s_j.
\]

Since $\mathrm{ord}_{i'}(\pi^*s_j) = \sum_i \mathrm{ord}_i(s_j) m_{ii'}$,
we should compare
\[
\llcorner t'_{i'} + \sum_i \mathrm{ord}_i(s_j) m_{ii'} \lrcorner
\]
and
\[
\min_t \{\sum_i (\llcorner t_i + \mathrm{ord}_i(s_j) \lrcorner)m_{ii'}
+ \llcorner t'_{i'} + \sum_i (1-t_i)m_{ii'} \lrcorner\}.
\]
We have
\[
\begin{split}
&\sum_i (\llcorner t_i + \mathrm{ord}_i(s_j) \lrcorner)m_{ii'}
+ \llcorner t'_{i'} + \sum_i (1-t_i)m_{ii'} \lrcorner
- \llcorner t'_{i'} + \sum_i \mathrm{ord}_i(s_j) m_{ii'} \lrcorner \\
&> \sum_i (\llcorner t_i + \mathrm{ord}_i(s_j) \lrcorner)m_{ii'}
+ \sum_i (1 - t_i - \mathrm{ord}_i(s_j) )m_{ii'} - 1
> - 1,
\end{split}
\]
hence
\[
\begin{split}
&\min_t \{\sum_i (\llcorner t_i + \mathrm{ord}_i(s_j) \lrcorner)m_{ii'}
+ \llcorner t'_{i'} + \sum_i (1-t_i)m_{ii'} \lrcorner\} \\
&\ge \llcorner t'_{i'} + \sum_i \mathrm{ord}_i(s_j) m_{ii'} \lrcorner.
\end{split}
\]
On the other hand, if we set 
$t_i = 1 - \mathrm{ord}_i(s_j) - \epsilon_i$ for $0 < \epsilon_i \ll 1$,
then
\[
\sum_i (\llcorner t_i + \mathrm{ord}_i(s_j) \lrcorner)m_{ii'}
+ \llcorner t'_{i'} + \sum_i (1-t_i)m_{ii'} \lrcorner\}
= \llcorner t'_{i'} + \sum_i \mathrm{ord}_i(s_j) m_{ii'} \lrcorner.
\]
Therefore, they are equal.
Since the minimum is attained at a value of $t$ which does not depend on $i'$
but only on $j$, we obtain the equality.
\end{proof}

Now we state the vanishing theorem for $\mathbb{Q}$-divisors in the relative setting.
This theorem will be used as an essential tool in the proof of the main theorem.

\begin{Thm}\label{vanishing}
In the situation of Definition~\ref{para}, let $L$ be a nef and big $\mathbb{Q}$-divisor 
on $X$ whose fractional part is 
supported on $B$.
Then
\[
H^p(X, \sum_t F^{t_1,\dots,t_h}(R^qf_*\omega_{Y/X}) 
\otimes \omega_X(\ulcorner L - \sum_i (1-t_i)B_i \urcorner)) = 0 
\]
for $p > 0$ and $q \ge 0$, where the sum is taken inside the sheaf 
$R^qf_*\omega_Y \otimes \mathcal{O}_X(\ulcorner L \urcorner)$.
\end{Thm}

\begin{proof}
By \cite{abel}, there exists a normal crossing divisor $\bar B$ such that
$B \le \bar B$ which satisfies the following:
there exists a finite surjective and Galois morphism $\pi: X' \to X$
from a smooth projective 
variety which is etale over $\bar X_0 = X \setminus \bar B$ 
and such that $\pi^*L$ has integral 
coefficients and all the local monodromies of 
$R^{d+q}f'_{0*}\mathbb{Z}_{Y'_0}$ are unipotent under the notation of Lemma~\ref{BC}.
We replace $B$ by $\bar B$ and
let $G$ be the Galois group of $\pi$.
By \cite{Kollar}~Theorem 2.1 and \cite{mK}~Theorem 3.3, we have 
\[
H^p(X', R^qf'_*\omega_{Y'} \otimes \mathcal{O}_{X'}(\pi^*L)) = 0 
\]
for $p > 0$.
Since $\pi^*(K_X + B) = K_{X'} + B'$,
we have by Lemma~\ref{BC} 
\[
\begin{split}
H^p(&X', \sum_t \pi^*F^{t_1,\dots,t_h}(R^qf_*\omega_{Y/X}) \\
&\otimes \mathcal{O}_{X'}(- \llcorner \sum_{i'} B'_{i'} 
+ \sum_i (1-t_i)\pi^*B_i \lrcorner + \pi^*(K_X + B + L))) = 0.
\end{split}
\]

We want to calculate the $G$-invariant part of the locally free sheaf
\[
\begin{split}
\mathcal{G} = F^{t_1,\dots,t_h}&(R^qf_*\omega_{Y/X}) \\
&\otimes \pi_*(\mathcal{O}_{X'}(- \llcorner \sum_{i'} B'_{i'} 
+ \sum_i (1-t_i)\pi^*B_i \lrcorner + \pi^*(K_X + B + L))).
\end{split}
\]
For this purpose, let $A$ be the largest divisor on $X$ such that
\[
\pi^*A \le \ulcorner - \sum_{i'} B'_{i'} - 
\sum_i (1-t_i)\pi^*B_i \urcorner + \pi^*(B + L). 
\]
This is equivalent to the condition
\[
\pi^*A < \pi^*(L - \sum_i (1-t_i)B_i) + \pi^*B.
\]
Hence we obtain
\[
A = \ulcorner L - \sum_i (1-t_i)B_i \urcorner.
\]
Since 
\[
0 = H^p(X, \mathcal{G})^G =  H^p(X, \sum_t F^{t_1,\dots,t_h}(R^qf_*\omega_{Y/X}) 
\otimes \omega_X(A)),
\]
our assertion is proved.
\end{proof}

\begin{Rem}
We note that the sum
\[
\sum_t F^{t_1,\dots,t_h}(R^qf_*\omega_{Y/X}) 
\otimes \omega_X(\ulcorner L - \sum_i (1-t_i)B_i \urcorner)
\]
is a locally
free sheaf because it is the $G$-invariant part of a locally free sheaf
as shown in the above proof, though the expression looks complicated.
It coincides with the subsheaf of $L^2$-sections of the locally free sheaf
$R^qf_*\omega_{Y/X} \otimes \omega_X(\ulcorner L \urcorner)$.
We can consider the non-vanishing problem for this sheaf.
\end{Rem}

\section{Proof of Theorem~\ref{main}}

Let $\{s_1, \cdots, s_k\}$ be the basis of 
$R^qf_*\omega_{Y/X}$ on a neighborhood $U$ of $x$
which is obtained in Lemma~\ref{basis}.
We shall prove that the image of the homomorphism
\[
H^0(X, R^qf_*\omega_Y \otimes L) \to 
R^qf_*\omega_Y \otimes L \otimes \kappa(x)
\]
contains $s_j \otimes \omega_X \otimes L \otimes \kappa(x)$ for any $j$.
Let us consider $\mathrm{ord}(s_j)$ as an effective $\mathbb{Q}$-divisor on $X$ by
setting the coefficients of the irreducible components of $B$ which do not intersect $U$
to be $0$.
By the assumption of the theorem, 
there exists an effective $\mathbb{Q}$-divisor $D$ such that
$D \sim_{\mathbb{Q}} \lambda L$ with $0 < \lambda < 1$,
$(X, \mathrm{ord}(s_j) + D)$ is properly log canonical at $x$, 
and that $\{x\}$ is a minimal log canonical center. 
By the perturbation of $D$, 
we may assume that $\{x\}$ is the only log canonical center which contains $x$,
and there exists only one log canonical place $E$ above the center $\{x\}$.
Let $\mu: X' \to X$ be a birational morphism from a smooth projective variety
such that $E$ appears as a smooth divisor on $X'$.
We write 
\[
\mu^*(K_X + \mathrm{ord}(s_j) + D) = K_{X'} + E + F,
\]
where the coefficients of $F$ are less than $1$.
We may assume that the union of the exceptional locus of $\mu$ and 
the support of $\mu^{-1}(B + D)$ is a normal crossing divisor.
Let $B' = \mu^*B_{\mathrm{red}} = \sum_{i'} B'_{i'}$. 
Since
\[
K_{X'} + (1 - \lambda)\mu^*L 
= \mu^*(K_X + L) - E - F + \mu^*\mathrm{ord}(s_j),
\]
we obtain by Theorem~\ref{vanishing}
\[
\begin{split}
H^1(X', &\sum_{t'} F^{t'_1,\dots,t'_{h'}}(R^qf'_*\omega_{Y'/X'}) 
\otimes \mathcal{O}_{X'}(\mu^*(K_X + L) \\
&- E + \ulcorner - F + \mu^*\mathrm{ord}(s_j) 
- \sum_{i'} (1 - t'_{i'})B'_{i'} \urcorner)) = 0.
\end{split}
\]
Since 
\[
\begin{split}
&\sum_{t'} F^{t'_1,\dots,t'_{h'}}(R^qf'_*\omega_{Y'/X'})  
\otimes \mathcal{O}_{X'}(\mu^*(K_X + L) \\
&+ \ulcorner - F + 
\mu^*\mathrm{ord}(s_j) - \sum_{i'} (1 - t'_{i'})B'_{i'} \urcorner))
\end{split}
\]
is a locally free sheaf on $X'$, we have a surjective homomorphism
\[
\begin{split}
H^0(X', &\sum_{t'} F^{t'_1,\dots,t'_{h'}}(R^qf'_*\omega_{Y'/X'}) 
\otimes \mathcal{O}_{X'}(\mu^*(K_X + L) \\
+ &\ulcorner - F + \mu^*\mathrm{ord}(s_j)
- \sum_{i'} (1 - t'_{i'})B'_{i'} \urcorner)) \\
\to H^0(&E, \sum_{t'} F^{t'_1,\dots,t'_{h'}}(R^qf'_*\omega_{Y'/X'}) 
\otimes \mathcal{O}_E(\mu^*(K_X + L) \\
&+ \ulcorner - F + \mu^*\mathrm{ord}(s_j)
- \sum_{i'} (1 - t'_{i'})B'_{i'} \urcorner)).
\end{split}
\]
We have 
\[
\mu_*(\ulcorner - F + \mu^*\mathrm{ord}(s_j)
- \sum_{i'} (1 - t'_{i'})B'_{i'} \urcorner) \le 0
\] 
if $0 \le t'_{i'} < 1$.
Hence 
\[
\begin{split}
&H^0(X', \sum_{t'} F^{t'_{1},\dots,t'_{h'}}(R^qf'_*\omega_{Y'/X'}) 
\otimes \mathcal{O}_{X'}(\mu^*(K_X + L) \\
&\quad + \ulcorner - F +  \mu^*\mathrm{ord}(s_j) 
- \sum_{i'} (1 - t'_{i'})B'_{i'} \urcorner)) \\
&\subset
H^0(X, R^qf_*\omega_{Y/X} \otimes \mathcal{O}_X(K_X + L)).
\end{split}
\]
We note that the $t'_{i'}$ need not be contained in the interval
$[0,1)$ in the above sum. 
On the other hand, if we define the $t'_{i'}$ by
\[
\sum_{i'} t'_{i'}B'_{i'} = \sum_{i'}(1 - \epsilon'_{i'})B'_{i'} - 
\mu^*\mathrm{ord}(s_j)
\]
for sufficiently small and positive numbers $\epsilon'_{i'}$, then
the divisor
\[
\ulcorner - F + \mu^*\mathrm{ord}(s_j) 
- \sum_{i'} (1 - t'_{i'})B'_{i'} \urcorner
\]
is effective and its support does not contain $E$, even if $E$ is contained
in $B'$.
Since $\mathrm{ord}(\mu^*s_j) = \mu^*\mathrm{ord}(s_j)$, we have
$\mu^*s_j \in F^{t'_1,\dots,t'_{h'}}(R^qf'_*\omega_{Y'/X'})$ for such
$t'_{i'}$. 
Hence
\[
\begin{split}
\mu^*s_j \otimes \mathcal{O}_E&(\mu^*(K_X + L)) \\
\in H^0(&E, \sum_{t'} F^{t'_1,\dots,t'_{h'}}(R^qf'_*\omega_{Y'/X'}) 
\otimes \mathcal{O}_E(\mu^*(K_X + L) \\
&+ \ulcorner - F + \mu^*\mathrm{ord}(s_j)
- \sum_{i'} (1 - t'_{i'})B'_{i'} \urcorner)).
\end{split}
\]
Therefore, $s_j \otimes \omega_X \otimes L \otimes \kappa(x)$ is 
contained in the image of $H^0(X, R^qf_*\omega_Y \otimes L)$.

Department of Mathematical Sciences, University of Tokyo, 

Komaba, Meguro, Tokyo, 153-8914, Japan 

kawamata@ms.u-tokyo.ac.jp

\end{document}